\input amstex
\documentstyle{amsppt}
\def \tagref#1{{\rm (#1)}}

\def \thmref#1{#1}
\def \dim#1{\text{dim}\,#1}
\font\oldsmallcaps=cmcsc10
\def\maple{{\oldsmallcaps maple\ }}
\NoBlackBoxes
\topmatter
\leftheadtext{A.O.L. Atkin and F.G. Garvan}
\rightheadtext{Relations between the rank and the crank}
\title
Relations between the ranks and cranks 
of partitions
\endtitle
\author A.O.L. Atkin and F.G. Garvan\endauthor
\address 
Department of Mathematics,
University of Illinois at Chicago,
851 S. Morgan Street,
Chicago, Illinois 60607-7045
\endaddress
\address Department of Mathematics, University of Florida, Gainesville,
Florida 32611\endaddress  
\email frank\@math.ufl.edu\endemail

\abstract 
New identities and congruences involving the ranks and cranks of partitions 
are proved.
The proof depends on a new partial differential equation connecting
their generating functions.
\endabstract

\dedicatory
Dedicated to the memory of Robert A. Rankin
\enddedicatory

\thanks 
Research of the second author supported in part by the NSF under grant
number DMS-9870052.
\endthanks
\endtopmatter

\document

\head 1. Introduction
\endhead
Dyson \cite{D1}, \cite{D3, p.52} defined the rank of a partition as the largest 
part minus
the number of parts. 
Let $N(m,n)$ denote the number of partitions of $n$ with rank
$m$, then
$$
\sum_m {N}(m,n) = p(n),
\tag"(1.1)"
$$
the number of partitions of $n$; and
$$
N(m,n) = N(-m,n),
\tag"(1.2)"
$$
using the classical conjugacy of partitions.

Andrews and Garvan \cite{A-G}
defined the crank of a partition. It is the largest part if the partition
contains no ones, and is otherwise 
the number of parts larger than the number of ones minus the number of ones.
Let ${M}(m,n)$ denote the number of partitions of $n$ with crank
$m$, then
$$
\sum_m {M}(m,n) = p(n),
\tag"(1.3)"
$$
the number of partitions of $n$; and
$$
M(m,n) = M(-m,n), 
\tag"(1.4)"
$$
for $n > 1$. A direct combinatorial proof of \tagref{1.4}
was found recently by Berkovich and Garvan \cite{B-G}.

We now state the generating function $R(z,q)$ for the rank.
We have
$$
\align
R(z,q) &= \sum_{n\ge0} \sum_{m} N(m,n) z^m q^n 
\tag"(1.5)"\\
&= 1 + \sum_{n=1}^\infty \frac{q^{n^2}}{(zq)_n(z^{-1}q)_n},
\tag"(1.6)"
\endalign
$$
and 
$$
\sum_{n\ge0} N(m,n) q^n 
 =\frac{1}{(q)_\infty} 
          \sum_{n=1}^\infty (-1)^{n-1}q^{\tfrac{n}{2}(3n-1) + |m|n}(1-q^n).
\tag"(1.7)"
$$
Here we are using the notation
$$
\align
(a)_n &= (a;q)_n = (1-a)(1-aq)\cdots(1-aq^{n-1}), \tag"(1.8)"\\
(q)_\infty &= (a;q)_\infty = \lim_{n\to\infty} (a)_n,
\endalign
$$
where $|q|<1$.

Below we state the generating function $C(z,q)$ for the crank.
If we amend the definition of $M(m,n)$ for $n=1$, then
the generating function can be given as an infinite product.
Accordingly, throughout this paper we assume
$$
M(-1,1)=-1,\, M(1,0)=M(1,1)=1,\text{and $M(m,1)=0$ otherwise}.
\tag"(1.9)"
$$
Then we have
$$
\align
C(z,q) &= \sum_{n\ge0} \sum_{m} M(m,n) z^m q^n 
\tag"(1.10)"\\
	&= \prod_{n=1}^\infty \frac{(1-q^n)}{(1-zq^n)(1-z^{-1}q^n)},
\tag"(1.11)"
\endalign
$$
and
$$
\sum_{n\ge0} M(m,n) q^n 
=\frac{1}{(q)_\infty} 
\sum_{n=1}^\infty (-1)^{n-1}q^{\tfrac{n}{2}(n-1) + |m|n}(1-q^n).
\tag"(1.12)"
$$
Equation \tagref{1.11} follows from \cite{A-G\rm, eq.~(1.11)} 
and \cite{A-G\rm, Thm~1}.     Equation \tagref{1.12} then
follows from \cite{G1\rm, eq.~(7.20)}.

The main result of this paper is a 
fundamental
partial differential equation (PDE) connecting
the rank and crank generating functions. See Theorem \thmref{1.1} below.
Let
$$
\align
R^{*}(z,q) &:= \frac{R(z,q)}{(1-z)},
\tag"(1.13)"\\
C^{*}(z,q) &:= \frac{C(z,q)}{(1-z)}.
\tag"(1.14)"
\endalign
$$
Define the differential operators
$$
\delta_z = z\frac{\partial}{\partial z},
\qquad
\delta_q = q\frac{\partial}{\partial q}.
\tag"(1.15)"
$$
\proclaim{Theorem 1.1}     We have
$$
z (q)_\infty^2 \left[ C^*(z,q) \right]^3
= \left(3\delta_q + \frac{1}{2} \delta_z + 
\frac{1}{2} \delta_z^2\right) R^*(z,q),
\tag"(1.16)"
$$
and
$$
z (q)_\infty^2 \left[ C(z,q) \right]^3
=
\left( 3 (1-z)^2 \delta_q + \frac{1}{2}(1-z)^2 \delta_z^2 
- \frac{1}{2}(z^2-1) \delta_z + z\right) R(z,q).
\tag"(1.17)"
$$
\endproclaim
We prove equation \tagref{1.16} in Section 2.
The Rank-Crank PDE \tagref{1.17} follows easily from \tagref{1.16}
by using \tagref{1.13} and \tagref{1.14}.

Let $N(k,t,n)$ denote the number of partitions of $n$ with rank congruent
to $k$ modulo $t$.
Then for $t=5$ or $t=7$
$$
N(k,t,n) = \frac{1}{t} \, p(n), \qquad 0\le k \le t-1;
\tag"(1.18)"
$$
for all $n$ satisfying $24n\equiv 1 \pmod{t}$.
These combinatorial results immediately imply Ramanujan's
partition congruences
$$
\align
p(5n+4) &\equiv 0 \pmod{5},
\tag"(1.19)"\\
p(7n+5) &\equiv 0 \pmod{7}.
\tag"(1.20)"
\endalign
$$

Let $M(k,t,n)$ denote the number of partitions of $n$ with crank congruent
to $k$ modulo $t$.
Then for $t=5$, $t=7$, or $t=11$
$$
M(k,t,n) = \frac{1}{t} \, p(n), \qquad 0\le k \le t-1;
\tag"(1.21)"
$$
for all $n$ satisfying $24n\equiv 1 \pmod{t}$.
These combinatorial results again imply Ramanujan's
partition congruences mod $5$ and mod $7$, and in addition his congruence
$$
p(11n+6) \equiv 0 \pmod{11}.
\tag"(1.22)"
$$

There are many more rank identities. For example,
$$
N(1,5,5n+1) = N(2,5,5n+1),
\tag"(1.23)"
$$
and others for the moduli $5$, $7$, $8$, $9$, and $12$.
The results for $5$ and $7$ were all found by Dyson \cite{D1}, \cite{D3, p.53}
and proved by Atkin and Swinnerton-Dyer \cite{A-SD}. The
results for $8$, $9$, and $12$, were found by Lewis \cite{L1}
and subsequently proved by Lewis and Santa-Gadea in
a series of papers \cite{SG1},
\cite{SG2}, and \cite{L-SG}.

There are similar identities for the crank. For example,
$$
M(0,8,2n+1)+M(1,8,2n+1)=M(3,8,2n+1)+M(4,8,2n+1),
\tag"(1.24)"
$$
and others for the moduli $5$, $7$, $8$, $9$, $10$, and $11$.
These were proved in \cite{G1}, \cite{G2}, and \cite{G3}.

There are identities between the rank and the crank. For example,
$$
M(4,9,3n) = N(4,9,3n),
\tag"(1.25)"
$$
and others for the moduli $5$, $7$, $8$, and $9$.
These results were proved in \cite{G2}, \cite{G3}, \cite{L1}, \cite{L3},
\cite{L-SG1}, and \cite{SG1}.

In this paper we consider 
linear relations modulo a prime $p$.
There are congruences for the rank to the moduli $11$ and $13$. 
For example,
$$
2 N(2,11,11n) + N(3,11,11n) 
 + 7N(4,11,11n) + N(5,11,11n)
 \equiv 0 \pmod{11}.
\tag"(1.26)"
$$
Results of this type are due to Atkin and Hussain \cite{A-H}
and O'Brien \cite{OB}.

It is a surprising fact that there is one
analogous relation for the 
crank for {\it every} prime $p$. This was a mystery to us until we
realized that these crank congruences follow
from the identity
$$
\sum_{k=1}^n k^2\, M(k,n) = n\,p(n),
\tag"(1.27)"
$$
due to Dyson \cite{D2}, who gave a combinatorial proof.

There is an extra linear congruence for the crank modulo $p$
for $p=41$, $53$, $83$, and $120667369$.

For each prime $p>13$ there are seven 
congruences involving both the rank and the crank  
modulo $p$. For example,
$$
\align
\qquad\qquad&6 N(0, 29, 29 n + 23) + 17 N(1, 29, 29 n + 23) + 24 N(2, 29, 29 n + 23)\tag"(1.28)"\\
&+ 18 N(3, 29, 29 n + 23) + 17 N(4, 29, 29 n + 23) + 14 N(5, 29, 29 n + 23) \\
&+ 22 N(6, 29, 29 n + 23) + 24 N(7, 29, 29 n + 23) + 2 N(9, 29, 29 n + 23) \\
&+15 N(10, 29, 29 n + 23) + 19 N(11, 29, 29 n + 23) + 18 N(12, 29, 29 n + 23)\\
& + 20 N(13, 29, 29 n + 23) + 16 N(14, 29, 29 n + 23) \\
&\equiv11M(0, 29, 29 n + 23)+17 M(1, 29, 29 n + 23) + 28 M(2, 29, 29 n + 23) \\
&\quad+ 26 M(4, 29, 29 n + 23) 
+ 6 M(5, 29, 29 n + 23) + 28 M(8, 29, 29 n + 23)\\
&\qquad\qquad\qquad\qquad\qquad\qquad\qquad\qquad\qquad\qquad \pmod{29}.
\endalign
$$
These congruences come from certain exact relations between the rank
and the crank.

For even $j\ge2$, we define
$$
\align
N_j(n) &= \sum_{k} k^j \, N(k,n),
\tag"(1.29)"\\
M_j(n) &= \sum_{k} k^j \, M(k,n).
\tag"(1.30)"
\endalign
$$
The following is the simplest exact relation
$$
N_4(n) 
= - (2n + \tfrac{2}{3})\, M_2(n) + \tfrac{8}{3}M_4(n) + (1 - 12n)\,N_2(n).
\tag"(1.31)"
$$
In fact, there are polynomials
$P_k(n)$ of degree $k-1$ and $Q_{k,j}(n)$ of degree $k-j$ (for $1\le j \le k$)
such that
$$
N_{2k}(n) = P_k(n) \, N_2(n) + \sum_{j=1}^k Q_{k,j}(n)\, M_{2j}(n),
\tag"(1.32)"
$$
for $k=2$, $3$, $4$, and $5$.
For $k=6$ there is no such relation.
For $k=7$ there is a similar relation but with an extra term $N_{12}(n)$.
See Theorem \thmref{5.1} below.
The proof of these exact relations depends on the Rank-Crank PDE \tagref{1.17}.

We call the functions $N_j$ and $M_j$ (defined in \tagref{1.29},
\tagref{1.30}) rank and crank moments, respectively.
We define the following generating functions
$$
\align
R_j &= R_j(q) = \sum_{n\ge1} N_j(n) q^n,
\tag"(1.33)"\\
C_j &= C_j(q) = \sum_{n\ge1} M_j(n) q^n,
\tag"(1.34)"
\endalign
$$
for even $j$. We find that
$$
\left.\delta_z^j R(z,q)\right|_{z=1} = 
\cases R_j, & \text{$j$ even},\\
0, & \text{$j$ odd},
\endcases
\tag"(1.35)"
$$
using \tagref{1.2} and \tagref{1.5}.
Similarly we find that
$$
\left.\delta_z^j C(z,q)\right|_{z=1} = 
\cases C_j, & \text{$j$ even},\\
0, & \text{$j$ odd}.
\endcases
\tag"(1.36)"
$$

In Section 2 we show how the Rank-Crank PDE \tagref{1.17}
follows from a certain elliptic-function identity \tagref{2.3}.
In Section 3 we prove some results for Eisenstein series, modular forms
and quasi-modular forms. In Section 4 we show how the crank moment functions
can be written in terms of Eisenstein series, and we derive some results for
the derivatives of crank moment functions. In Section 5 we show how the
Rank-Crank PDE and certain results for the derivatives of Eisenstein
series lead to exact relations between rank and crank moments.
As a bonus we show that the $23$rd power of the Dedekind eta-function
can be written in terms of rank and crank moments.
In Section 6 we consider congruence relations between rank and crank moments.

\vfill
\eject

\head 2. Proof of the rank-crank PDE
\endhead

The rank-crank PDE follows easily from an identity in \cite{A-SD}.
Define
$$
J(z,q) := \prod_{n=1}^\infty (1 - z^{-1}q^n)(1 - z q^{n-1}),
\tag"(2.1)"
$$
and
$$
S(z,\zeta,q) := \sum_{n=-\infty}^\infty
\frac{ (-1)^n \zeta^n q^{3n(n+1)/2}}{1-z q^n}.
\tag"(2.2)"
$$
Then
$$
\zeta^3 S(z\zeta,\zeta^3,q) + S(z\zeta^{-1},\zeta^{-3},q)
 - \zeta \frac{ J(\zeta^2,q)}{J(\zeta,q)} S(z,1,q)
=
\frac{ J(\zeta,q) J(\zeta^2,q) (q)_\infty^2}
{J(\zeta z,q) J(z,q) J(z\zeta^{-1},q)}.
\tag"(2.3)"
$$
This identity is equation (5.1) in \cite{A-SD\rm, p.94} and was one of the key 
identities
required to prove Dyson's results for the rank modulo $5$ and $7$.
We let $g(\zeta)$ denote either side of \tagref{2.3} as a function
of $\zeta$. By considering the right side of \tagref{2.3} we see
that $g(\zeta)$ has a double zero at $\zeta=1$
and that
$$
g''(1) = 4 (q)_\infty^3 \left[ C^*(z,q) \right]^3,
\tag"(2.4)"
$$
where $C^*(z,q)$ is defined in \tagref{1.14}.
We let $h(\zeta)$ be the sum of the first two terms on the left side
of \tagref{2.3}; i.e.,
$$
h(\zeta) = 
\zeta^3 S(z\zeta,\zeta^3,q) + S(z\zeta^{-1},\zeta^{-3},q).
\tag"(2.5)"
$$
We find that
$$
\align
h''(1)
&= \sum_{n=-\infty}^\infty (-1)^n q^{3n(n+1)/2}
\left(
\frac{ 6(3n^2+n+1) }{1-zq^n}
+\frac{ 4 z (3n+1) q^n}{(1-zq^n)^2}
+\frac{ 4 z^2 q^{2n}}{(1-zq^n)^3}
\right)\\
&=(6 + 12 \delta_q + 6 \delta_z + 2\delta_z^2) S(z,1,q).
\tag"(2.6)"
\endalign
$$

We let $j(\zeta)$ be the third term on the right side of \tagref{2.3}
$$
j(\zeta) = 
 - \zeta \frac{ J(\zeta^2,q)}{J(\zeta,q)} S(z,1,q).
\tag"(2.7)"
$$
We find that
$$
j''(1) = -2 \left(
1 - 6 \sum_{n\ge1} \frac{q^n}{(1-q^n)^2}\right)
S(z,1,q) = -2 (1 - 6\Phi_1(q)) S(z,1,q),
\tag"(2.8)"
$$
where
$$
\Phi_1(q) = \sum_{n=1}^\infty \frac{n q^n}{1-q^n}.
\tag"(2.9)"
$$
The functions $\Phi_j$ are defined below in \tagref{3.1}.
We define
$$
P(q) = \sum_{n=0}^\infty p(n) q^n = \frac{1}{(q)_\infty}.
\tag"(2.10)"
$$
By differentiating logarithmically with respect to $q$ we obtain the
well-known identity
$$
\delta_q P(q) = \Phi_1(q) \, P(q),
\tag"(2.11)"
$$
or
$$
\delta_q (q)_\infty = - \Phi_1(q) (q)_\infty.
\tag"(2.12)"
$$
The following identity is equation (7.10) in \cite{G1}:
$$
z S(z,1,q) = (q)_\infty(-1 + R^*(z,q)).
\tag"(2.13)"
$$
We apply $\delta_q$ to both sides of \tagref{2.13} and use \tagref{2.12}
to find that
$$
z\delta_q S(z,1,q) = (q)_\infty \delta_q R^*(z,q) - z\Phi_1(q) S(z,1,q).
\tag"(2.14)"
$$
Similarly we find that
$$
\align
z\delta_z S(z,1,q) &= (q)_\infty \delta_z R^*(z,q) - z S(z,1,q),\tag"(2.15)"\\
z\delta_z^2 S(z,1,q) &= (q)_\infty \left(\delta_z^2 - 2\delta_z\right) R^*(z,q)
+ z S(z,1,q).\tag"(2.16)"
\endalign
$$
Now 
$$
g''(1) = h''(1) + j''(1).
\tag"(2.17)"
$$
Using \tagref{2.4}, \tagref{2.6}, \tagref{2.8}, and \tagref{2.13}--\tagref{2.16}
this equation becomes
$$
z (q)_\infty^2 \left[ C^*(z,q) \right]^3
= \left(3\delta_q + \frac{1}{2} \delta_z 
+ \frac{1}{2} \delta_z^2\right) R^*(z,q).
\tag"(2.18)"
$$
 From \tagref{1.13} and \tagref{1.14} 
we  find that
$$
\align
\delta_z R^*(z,q) &= \frac{ \delta_z R(z,q) + z R^*(z,q)}{1-z},
\tag"(2.19)"\\
\delta_z^2 R^*(z,q) &= \frac{ \delta_z^2 R(z,q) 
                       + 2z \delta_z R^*(z,q) + z R^*(z,q)}{1-z},
\tag"(2.20)"\\
\delta_q R^*(z,q) &= \frac{\delta_q R(z,q)}{1-z}.
\tag"(2.21)"
\endalign
$$
Using (2.19)--(2.21) 
we can write \tagref{2.18} in terms of $C(z,q)$ and $R(z,q)$:
$$
z (q)_\infty^2 \left[ C(z,q) \right]^3
=
\left( 3 (1-z)^2 \delta_q + \frac{1}{2}(1-z)^2 \delta_z^2 
- \frac{1}{2}(z^2-1) \delta_z + z\right) R(z,q),
\tag"(2.22)"
$$
which is the rank-crank PDE.

\vfill
\eject

\head 3. Eisenstein series, modular forms and derivatives
\endhead

Following Ramanujan \cite{Ram\rm, p.163} we define
$$
\Phi_j(q) = \sum_{n=1}^\infty \frac{n^j q^n}{1-q^n} = \sum_{m,n\ge1} n^j q^{nm}
= \sum_{n=1}^\infty \sigma_j(n) q^n,
\tag"(3.1)"
$$
for $j\ge1$ odd and where $\sigma_j(n) = \sum_{d\mid n} d^j$.
As usual we let $q=\exp(2\pi i\tau)$ where $\tau$ is in the complex
upper half-plane $\Cal H$ so that $|q|<1$.
For $n$ even the Eisenstein series $E_n(\tau)$ is
defined by
$$
\align
E_n(\tau) &= 1 - \frac{2n}{B_n} \sum_{k=1}^\infty \sigma_{n-1}(k) e^{2\pi ik\tau},\\
&\tag"(3.2)"\\
&= 1 - \frac{2n}{B_n} \Phi_{n-1}(q),
\endalign
$$
where the Bernoulli numbers $B_n$ are defined by
$$
\frac{x}{e^x-1} = \sum_{n=0}^\infty B_n \frac{x^n}{n!}.
$$
Ramanujan \cite{R\rm, p.140} considered in particular the Eisenstein series
$$
\align
E_2 &= 1 - 24 \, \Phi_1,\\
E_4 &= 1 + 240\, \Phi_3,\tag"(3.3)"\\
E_6 &= 1 - 540\, \Phi_5.  
\endalign
$$
For even $n\ge 4$, $E_n$ is
a modular form of weight $n$ for the full modular group 
$\Gamma=\text{SL}_2(\Bbb Z)$ \cite{Ran}.
$E_2$ is not a modular form, but is transformed by the generators
of $\Gamma$ according to
$$
\align
E_2(\tau + 1) &= E_2(\tau),\tag"(3.4)"\\
\tau^{-2} E_2(-1/\tau) &= E_2(\tau) + \frac{12}{2\pi i\tau}.
\endalign
$$
See \cite{K\rm, p.113}.
Ramanujan \cite{Ram\rm, p.165} found that
$$
\align
\delta_q( E_2 ) &= \frac{E_2^2 - E_4}{12},\\
\delta_q( E_4 ) &= \frac{E_2 E_4 - E_6}{3}, \tag"(3.5)"\\
\delta_q( E_6 ) &= \frac{E_2 E_6 - E_4^2}{2},
\endalign
$$
where
$$
\delta_q = q\,\frac{d}{dq}.
$$
More generally, it is known that if $f$ is a modular form of
weight $k$ then
$$
12 \delta_q(f) - k E_2 f
$$
is a modular form of weight $(k+2)$. See \cite{SD\rm, p.19}.

Following Serre \cite{S\rm, p.88}, we let $\Cal M_{k}$ denote the vector space of
modular forms of weight $2k$. Then
$$
\text{dim} \Cal M_{k} = 
\cases
[k/6], & \text{if $k\equiv 1\pmod{6}$},\\
[k/6]+1, & \text{otherwise},
\endcases
\tag"(3.6)"
$$
and the
set
$$
\{ E_4^{a} E_6^{b}\,:\, \text{$2a+3b=k$ with  $a$ and $b$ nonnegative 
integers}\}
\tag"(3.7)"
$$
forms a basis for $\Cal M_{k}$. See Serre \cite{S\rm, p.88}.
Thus any Eisenstein series $E_{2n}$ (for $n\ge2$)
can be written in terms of $E_4$ and $E_6$. This can be done explicitly using
well-known recurrences. See \cite{Ap\rm, pp.12-13} and \cite{B\rm, pp.331-332}.

The graded algebra of modular forms is given by
$$
\Cal M = \sum_{k=0}^\infty \Cal M_k = \Bbb C[E_4,E_6].
\tag"(3.8)"
$$
See Serre \cite{S,\rm p.89}.
We need to extend this algebra to include $E_2$. We say that $f$ is a
{\it quasi-modular form} if it is in the algebra generated
by $E_2$ and $\Cal M$. We extend the definition of weight by defining the
weight of $E_2$ to be $2$. 
Let $n$ be a nonnegative integer. 
Then      
the space of quasi-modular forms of weight $\le 2n$ is 
$$
\Cal E_n = \left\{ \sum_{j=0}^n f_j E_2^j \,:\, 
            f_j \in \sum_{k=0}^{n-j} \Cal M_k\right\},
\tag"(3.9)"
$$
which is clearly a vector space over $\Bbb C$.
Below in Corollary \thmref{3.6} we give a basis for $\Cal E_n$.

Quasi-modular forms were first studied systematically by Kaneko and 
Zagier \cite{K-Z}.
We need some independence results for modular forms and 
quasi-modular forms. 
The results for quasi-modular forms 
follow from \cite{K-Z; Proposition 1(b), p.167}. We have included proofs
for completeness since they are elementary and the details of the
relevant proof in \cite{K-Z} are omitted.               
\proclaim{Lemma 3.1} Let $n$ be a nonnegative integer. 
Suppose that
$$
f_k\,:\, \Cal H \longrightarrow \Bbb C, \text{ where } f_k(\tau + 1) = f_k(\tau), \qquad (0 \le k \le n)
\tag"(3.10)"
$$
for all $\tau \in \Cal H$, and 
$$
\sum_{k=0}^n \tau^k f_k(\tau) = 0,
\tag"(3.11)"
$$
for all $\tau \in \Cal H$. Then
$$
f_0(\tau)=f_1(\tau)= \cdots = f_n(\tau) = 0,
\tag"(3.12)"
$$
for all $\tau \in \Cal H$. 
\endproclaim
\demo{Proof} 
Suppose \tagref{3.10} and \tagref{3.11} hold. Let $\tau\in\Cal H$ be
arbitrary but fixed. Then
$$
\sum_{k=0}^n (\tau+m)^k f_k(\tau) = 0,
\tag"(3.13)"
$$
for all integers $m$. Hence the polynomial
$$
p(z) = \sum_{k=0}^n f_k(\tau) z^k
$$
has infinitely many zeros. The result follows. \qed
\enddemo

\proclaim{Corollary 3.2} Non-zero modular forms of different weights are
linearly independent over $\Bbb C$.
\endproclaim
\demo{Proof}
Suppose that there are complex constants $c_k$ such that
$$
\sum_{k=0}^n c_k f_k(\tau) = 0,
\tag"(3.14)"
$$
for all $\tau \in \Cal H$, 
where
$f_k(\tau)$ is a modular form of weight $k$.
We apply $\tau \to -1/\tau$ to obtain
$$
\sum_{k=0}^n c_k \tau^k f_k(\tau) = 0,
$$
for all $\tau \in \Cal H$. The functions $f_k$ satisfy \tagref{3.10} so by 
Lemma \thmref{3.1} we have
$$
c_0 f_0(\tau)=c_1 f_1(\tau)= \cdots = c_n f_n(\tau) = 0,
$$
and the result follows. \qed
\enddemo

\proclaim{Proposition 3.3}
Let $n$ be a nonnegative integer. 
Suppose that
$$
f(\tau) := \sum_{j=0}^n f_j(\tau) E_2^j(\tau) = 0,
\tag"(3.15)"
$$
for all $\tau \in \Cal H$, and $f$ is a quasi-modular form of
weight $\le 2n$, so that each $f_j$ is a sum of modular forms
of weight $\le 2n-2j$. Then
$$
f_0=f_1= \cdots = f_n = 0.
$$
\endproclaim
\demo{Proof}
Suppose \tagref{3.15} holds and suppose for $0\le j \le n$
$$
f_j = f_{j,0} + f_{j,4} + \cdots + f_{j,2n-2j},
$$
where each $f_{j,k}$ is a modular form of weight $k$. Hence
$$
\sum_{j=0}^n\left(
f_{j,0} + \sum_{k=2}^{n-j} f_{j,2k} \right) E_2^j = 0.
\tag"(3.16)"
$$
Applying $\tau\to -1/\tau$ we obtain
$$
\sum_{j=0}^n\left(
f_{j,0} + \sum_{k=2}^{n-j} \tau^{2k} f_{j,2k} \right) 
\left(\tau^2 E_2 + \alpha \tau\right)^j = 0,
\tag"(3.17)"
$$
where
$$
\alpha = - \frac{6i}{\pi},
\tag"(3.18)"
$$
by \tagref{3.4}.
We rewrite \tagref{3.17} as
$$
\sum_{j=0}^n\left(
f_{j,0} + \sum_{k=2}^{n-j} \tau^{2k} f_{j,2k} \right) 
\tau^j
\sum_{\ell=0}^j 
\binom{j}{\ell} \tau^{\ell} E_2^{\ell} \alpha^{j-\ell}
= 0,
$$
whence
$$
\align
&f_{0,0} + \alpha f_{1,0} \tau \tag"(3.19)"\\
&+ \sum_{m=2}^{2n} \left(
\sum\Sb 0\le j \le n \\ 0\le \ell \le j \\ j+\ell = m\endSb
f_{j,0} \binom{j}{\ell} E_2^{\ell} \alpha^{j-\ell} 
 +
\sum\Sb 0\le j \le n \\ 2 \le k \le n-j\\ 0\le \ell \le j \\ j+\ell+2k = m\endSb
\binom{j}{\ell}\alpha^{j-\ell} f_{j,2k} E_2^{\ell} \right) \tau^m = 0.
\endalign
$$
Since the $f_{j,k}$ and $E_2$ satisfy \tagref{3.10}, by Lemma \thmref{3.1}  we 
have
$$
\align
&f_{0,0} = 0,\\
&f_{1,0} = 0,\tag"(3.20)"\\
&
\sum\Sb 0\le j \le n \\ 0\le \ell \le j \\ j+\ell = m\endSb
f_{j,0} \binom{j}{\ell} E_2^{\ell} \alpha^{j-\ell} 
 +
\sum\Sb 0\le j \le n \\ 2 \le k \le n-j\\ 0\le \ell \le j \\ j+\ell+2k = m\endSb
\binom{j}{\ell}\alpha^{j-\ell} f_{j,2k} E_2^{\ell} = 0,
\quad 2 \le m \le 2n.
\endalign
$$
There are $2n+1$ equations in \tagref{3.20}.
Each $f_{j,2k}$ occurs in some equation. We prove that
$$
f_{j,2k}=0,
\tag"(3.21)"
$$
for all $j$, $k$ such that $f_{j,2k}$ occurs in the $m$-th equation 
of \tagref{3.20}.
We proceed by induction on $m$. The result holds for $m=0$, and $m=1$.
Suppose it holds for $m < M$ where $M$ is a fixed positive integer $\le 2n$.
We consider the $M$-th equation. 
It is clear that if a term $f_{j,2k} E_2^{\ell}$ with $\ell \ge 1$
occurs in this equation, then $f_{j,2k}$ must have appeared in a previous
equation, and so is zero by the induction hypothesis. The remaining terms
in the $M$-th equation all have $\ell=0$, and are true modular forms, and it
is easy to see that they have different weights. Hence they are all
zero by Corollary \thmref{3.2}, 
and thus the result holds for $m=M$, and so for all
$m$ with $0\le m\le2n$ by induction. \qed
\enddemo

\proclaim{Corollary 3.4} Let $n$ be a nonnegative integer.
The set 
$$
\{
E_2^a E_4^b E_6^c \,:\, \text{$a+2b+3c\le n$ with $a$, $b$, $c$ nonnegative
integers}\}
\tag"(3.22)"
$$
is a basis for the space of quasi-modular forms of weight $\le 2n$.
\endproclaim

In the next section we shall see that the crank moment functions
$C_a$ can be written in terms of the $\Phi_j$ but are not modular forms.
Thus we need some results for the $\Phi_j$.
Using \tagref{3.3} and \tagref{3.5} we find that
$$
\align
\delta_q(\Phi_1) &=
\frac{5}{6}\,\Phi_{{3}}+\frac{1}{6}\,\Phi_{{1}}-2\,{\Phi_{{1}}}^{2},\\
\delta_q(\Phi_3) &=
{\frac {7}{10}}\,\Phi_{{5}}+\frac{1}{3}\,\Phi_{{3}}-\frac{1}{30}\,\Phi_{{1}}-8\,\Phi_ {{1}}\Phi_{{3}},\tag"(3.23)"\\
\delta_q(\Phi_5) &=
\frac{1}{2}\,\Phi_{{5}}+\frac{1}{42}\,\Phi_{{1}}-12\,\Phi_{{1}}\Phi_{{5}}+{\frac {10} {21}}\,\Phi_{{3}}+{\frac {400}{7}}\,{\Phi_{{3}}}^{2}.
\endalign
$$

We define the weight $\omega$ of a monomial $\Phi_1^a \Phi_3^b \Phi_5^c$
as
$$
\omega(\Phi_1^a \Phi_3^b \Phi_5^c) = 2a + 4b + 6c.
\tag"(3.24)"
$$
It is clear from \tagref{3.3} that the monomial is a sum of
quasi-modular forms of different weights, the form
$E_2^a E_4^b E_6^c$ having the highest weight $\omega$.
 From \tagref{3.2}, \tagref{3.3} and the fact the set in \tagref{3.7}
forms a basis for $\Cal M_n$, we see that for $n>1$
there are constants $\alpha_{b,c}$ such that
$$
\Phi_{2n-1} = \sum_{0<2b+3c\le n} \alpha_{b,c} \Phi_3^b \Phi_5^c.
\tag"(3.25)"
$$
In the sum $2b+3c$ is positive since $\Phi_{2n-1}(q)=0$ for $q=0$.
For example,
$$
\Phi_{11} = \frac{1}{13}\left(
63(
\Phi_{{3}}+240\,{\Phi_{{3}}}^{2}+19200\,{\Phi_{{3}}}^{3}+200\,{\Phi_{{5}}}^{2})
-50 \Phi_5\right).
\tag"(3.26)"
$$

For $n\ge1$ we let $\Cal V_n$ be the $\Bbb C$-vector space spanned by
monomials $\Phi_1^a \Phi_3^b \Phi_5^c$ with $a + 2b + 3c =n$.
We define
$$
\Cal W_n = \sum_{k=1}^n \Cal V_n;
\tag"(3.27)"
$$
i.e., $\Cal W_n$ is the $\Bbb C$-vector space spanned by
monomials $\Phi_1^a \Phi_3^b \Phi_5^c$ with $0< a + 2b + 3c \le n$.
Using \tagref{3.23} and induction we can show that
$$
\delta_q \left(\Cal W_n\right) \subset \Cal W_{n+1},
\tag"(3.28)"
$$
so that
$$
\delta_q^m \left(\Cal W_n\right) \subset \Cal W_{n+m},
\tag"(3.29)"
$$
for $m\ge 0$.

\proclaim{Theorem 3.5} Let $n$ be a nonnegative integer.
The set 
$$
\{
\Phi_1^a \Phi_3^b \Phi_5^c \,:\, 
\text{$0<a+2b+3c\le n$ with $a$, $b$, $c$ nonnegative
integers}\}
\tag"(3.30)"
$$
is a basis for $\Cal W_n$. 
\endproclaim
\demo{Proof}
The set in \tagref{3.30} spans $\Cal W_n$ by definition. We show that the
set is linearly independent, by showing a slightly larger set
is linearly independent. By equation \tagref{3.3} and Corollary \thmref{3.4}
the set of monomials $\Phi_1^a \Phi_3^b \Phi_5^c $, where
$0\le a+2b+3c\le n$, spans $\Cal E_n$ and the number of such monomials
is equal to $\dim{\Cal E_n}$.  Hence, these monomials form a 
basis for $\Cal E_n$ and are linearly independent. \qed
\enddemo
\proclaim{Corollary 3.6} For $n\ge1$
$$
\align
\text{\rm dim}\,\Cal V_n = \sum_{k=0}^n \text{\rm dim}\,\Cal M_k,\tag"(3.31)"\\
\text{\rm dim}\,\Cal W_n = n + \sum_{k=2}^n (n-k+1)\text{\rm dim}\,\Cal M_k.
\tag"(3.32)"
\endalign
$$
\endproclaim
\demo{Proof} 
The result follows from Theorem \thmref{3.5} and the fact that the set in 
\tagref{3.7} forms a basis for $\Cal M_k$.  \qed
\enddemo

We have the following table.
$$
\matrix
k & \text{dim}\,\Cal M_k & \text{dim}\,\Cal V_k & \text{dim}\,\Cal W_k\\
\noalign{\smallskip\hrule\medskip}
1&0&1&1\\ \noalign{\medskip}
2&1&2&3\\\noalign{\medskip}
3&1&3&6\\\noalign{\medskip}
4&1&4&10\\\noalign{\medskip}
5&1&5&15\\\noalign{\medskip}
6&2&7&22\\\noalign{\medskip}
7&1&8&30\\\noalign{\medskip}
8&2&10&40\\\noalign{\medskip}
9&2&12&52\\\noalign{\medskip}
10&2&14&66
\endmatrix
$$

\vfill
\eject

\head 4. Crank moments
\endhead

We now prove that the crank moment functions $C_a$ for $a$ even
can be written in terms of $P$ and the $\Phi_j$. From \tagref{1.11} we have
$$
\delta_z C(z,q) = L(z,q) C(z,q),
\tag"(4.1)"
$$
where
$$
\align
L(z,q) &= \sum_{n\ge1} \left(\frac{zq^n}{1-zq^n} - 
\frac{z^{-1}q^n}{1-z^{-1}q^n}\right),
\tag"(4.2)"\\
&= \sum_{m,n\ge1} \left(z^m q^{mn} - z^{-m} q^{mn}\right),
\endalign
$$
so that
$$
\delta_z^j L(z,q) = \sum_{m,n\ge1} \left(m^j z^m q^{mn} - 
                        (-m)^j z^{-m} q^{mn}\right),
\tag"(4.3)"
$$
and
$$
\left.\delta_z^j L(z,q)\right|_{z=1} = 
\cases 0, & \text{$j$ even},\\
2\Phi_j, & \text{$j$ odd}.
\endcases
\tag"(4.4)"
$$

Assume $a$ is even and apply $\delta_z^{a-1}$ to both sides of \tagref{4.1}.
Then
$$
\delta_z^a C = \sum_{j} \binom{a-1}{j} \delta_z^j (L)\, \delta_z^{a-1-j} (C).
\tag"(4.5)"
$$
Setting $z=1$ and using \tagref{4.4} and \tagref{1.36} we obtain the following
recurrence
$$
C_a = 2\sum_{j=1}^{\frac{a}{2} - 1}
      \binom{a-1}{2j-1}\,\Phi_{2j-1}\,C_{a-2j}
      + 2 \Phi_{a-1}\,P.
\tag"(4.6)"
$$  
We compute some examples.
$$
\align
C_2 &= 2\,P\, \Phi_{{1}},\tag"(4.7)"\\
C_4 &=2\,{ P}\,\left (\Phi_{{3}}+6\,{\Phi_{{1}}}^{2}\right ),\\
C_6 &=
2\,{ P}\,\left (\Phi_{{5}}+30\,\Phi_{{3}}\Phi_{{1}}+60\,{\Phi_{{1} }}^{3}\right ),
\\
C_8 &=
2\,{ P}\,\left (\Phi_{{7}}+56\,\Phi_{{5}}\Phi_{{1}}+840\,\Phi_{{3} }{\Phi_{{1}}}^{2}+840\,{\Phi_{{1}}}^{4}+70\,{\Phi_{{3}}}^{2}\right ).
\endalign
$$
Using induction and \tagref{4.6} we can show that for $n\ge1$
there are integers $\alpha_{a_1,a_2,\dots,\alpha_n}$ such that
$$
C_{2n} = 2\, P\, \sum_{a_1+2a_2 + \cdots + n a_n = n} \alpha_{a_1,a_2,\dots,\alpha_n} \Phi_1^{a_1} \Phi_3^{a_2}
\cdots \Phi_{2n-1}^{a_n}.
\tag"(4.8)"
\endalign
$$
By \tagref{2.11} and \tagref{4.7} we have
$$
C_2 = 2 \delta_q P,
\tag"(4.9)"
$$
or
$$
M_2(n) = \sum_k k^2\, M(k,n) = 2 n p(n).
\tag"(4.10)"
$$
A combinatorial proof of \tagref{4.10} was found by Dyson \cite{D2}.

We need the following
\proclaim{Lemma 4.1}
For $m\ge1$, there exists a $\Phi\in\Cal W_{m}$ such that
$$
\delta_q^m\left(P\right) = P \, \Phi.  
\tag"(4.11)"
$$
\endproclaim
\demo{Proof} We proceed by induction on $m$. From \tagref{2.11}
$$
\delta_q(P) = P \, \Phi_1,
\tag"(4.12)"
$$
and the result holds for $m=1$. Suppose the result holds for $m=a$; i.e., 
$$
\delta_q^a\left(P\right) = P \, \Phi,  
\tag"(4.13)"
$$
for some $\Phi\in\Cal W_{a}$. Then
$$
\align
\delta_q^{a+1}\left(P\right) &= P \, \delta_q(\Phi) + \delta_q(P)\, \Phi,  \\
&= P( \delta_q(\Phi) + \Phi_1\, \Phi),
\endalign
$$
and the result holds for $m=a+1$ since $\delta_q(\Phi)\in\Cal W_{a+1}$
by \tagref{3.29}. The result for general $m$ follows by induction. \qed
\enddemo

We may calculate $\delta_q^a(P)$ by using \tagref{4.12}
and the recurrence
$$
\delta_q^a(P) = \sum_{j=0}^{a-1} \binom{a-1}{j} \, \delta_q^j(\Phi_1)
\, \delta_q^{a-j-1}(P),
\tag"(4.14)"
$$
which is obtained by applying $\delta_q^{a-1}$ to both sides of \tagref{4.12}.
For example, we have
$$
\align
\delta_q^2(P) &=-\frac{1}{6}\,{P}\,\left (6\,{\Phi_{{1}}}^{2}-5\,\Phi_{{3}}-\Phi_{{1}} \right ),\tag"(4.15)"\\
\delta_q^3(P) &=
\frac{1}{12}\,{P}\,\left (36\,{\Phi_{{1}}}^{3}-90\,\Phi_{{1}}\Phi_{{3}}-6 \,{\Phi_{{1}}}^{2}+7\,\Phi_{{5}}+5\,\Phi_{{3}}\right ).
\endalign
$$

\proclaim{Theorem 4.2}
For $m\ge0$ and $n\ge1$, there exists $\Phi\in\Cal W_{n+m}$ such that
$$
\delta_q^m\left(C_{2n}\right) = P \, \Phi,
\tag"(4.16)"
$$
where $\Cal W_k$ is defined by \tagref{3.27}.
\endproclaim
\demo{Proof}
Let $n\ge1$. We proceed by induction on $m$. The result is true for $m=0$
using \tagref{4.8} and \tagref{3.25}. The remainder of the proof is analogous
to that of Lemma \thmref{4.1}. \qed
\enddemo

We give some examples.
$$
\align
\delta_q(C_2) &=
-\frac{1}{3}\,{P}\,\left (6\,{\Phi_{{1}}}^{2}-5\,\Phi_{{3}}-\Phi_{{1}}\right ),
\\
\delta_q^2(C_2) &=
\frac{1}{6}\,{P}\,\left (36\,{\Phi_{{1}}}^{3}-90\,\Phi_{{1}}\Phi_{{3}}-6\,{\Phi _{{1}}}^{2}+7\,\Phi_{{5}}+5\,\Phi_{{3}}\right ),\tag"(4.17)"
\\
\delta_q(C_4) &=
-\frac{1}{15}\,{P}\,\left (-90\,\Phi_{{1}}\Phi_{{3}}-21\,\Phi_{{5}}-10\,\Phi_{{ 3}}+\Phi_{{1}}+540\,{\Phi_{{1}}}^{3}-60\,{\Phi_{{1}}}^{2}\right ).
\endalign
$$

\head 5. Relations between rank and crank moments
\endhead

Let $a$ be even. After applying $\delta_z^a$ to both sides of the
rank-crank PDE \tagref{2.22}, setting $z=1$ and using \tagref{4.9} we find
that
$$
\align
&\sum_{i=0}^{a/2-1}\binom{a}{2i}
\sum\Sb \alpha+\beta+\gamma=a-2i\\ \text{$\alpha$, $\beta$, $\gamma$ even $\ge0$}\endSb
\binom{a-2i}{\alpha,\beta,\gamma}\,C_\alpha\,C_\beta\,C_\gamma\,P^{-2}
- 3\left(2^{a-1}-1\right) C_2 \\
&= \frac{1}{2}(a-1)(a-2) R_a
+ 6\sum_{i=1}^{a/2-1} \binom{a}{2i}\left(2^{2i-1}-1\right) \delta_q(R_{a-2i})
\tag"(5.1)"
\\
&\quad + \sum_{i=1}^{a/2-1}\left[
\binom{a}{2i+2}\left(2^{2i+1}-1\right) - 2^{2i}\binom{a}{2i+1} + \binom{a}{2i}
\right] R_{a-2i}.
\endalign
$$
For $a=2$ we obtain $0=0$. For $a=4$ we obtain
$$
C_{{4}}+6\,{\frac {{C_{{2}}}^{2}}{{P}}}-C_{{2}}=R _{{4}}-\,R_{{2}}
+12\,\delta_q(R_{{2}}).
\tag"(5.2)"
$$
Using \tagref{4.7}, \tagref{4.17} and \tagref{5.2} we find that
$$
-\frac{2}{3}\,C_{{2}}-2\,\delta_q(C_{{2}})+\frac{8}{3}\,C_{{4}}=R_{{4}}-R_{{2}}+12\,\delta_q({R}_{{2}}),
\tag"(5.3)"
$$
or 
$$
N_{4}(n)
=\frac{2}{3}\,\left (-3\,n-1\right )\,M_{2}(n)+\frac{8}{3}\,M_{4}( n)+\left (-12\,n+1\right )\,N_{2}(n),
\tag"(5.4)"
$$
for $n\ge0$.

Similarly, for $2$, $3$, $4$, and $5$, there are polynomials
$P_k(n)$ of degree $k-1$ and $Q_{k,j}(n)$ of degree $k-j$ (for $1\le j \le k$)
such that
$$
N_{2k}(n) = P_k(n) \, N_2(n) + \sum_{j=1}^k Q_{k,j}(n)\, M_{2j}(n),
\tag"(5.5)"
$$
for $n\ge0$. For $k=6$ there is no such relation. For $k=7$ there is 
a similar relation but with an extra term $N_{12}(n)$. These
relations are given in the following

\proclaim{Theorem 5.1}
For $n\ge0$ we have
$$
N_{4}(n)
=\frac{2}{3}\,\left (-3\,n-1\right )\,M_{2}(n)+\frac{8}{3}\,M_{4}( n)+\left (-12\,n+1\right )\,N_{2}(n),
\tag"(5.6)"
$$
$$
\align
N_{6}(n)&={\frac {2}{33}}\,\left (324\,{n}^{2}+69\,n-10\right )\,M_{2}(n)+{\frac {20}{33}}\,\left (-45\,n+4\right )\,M_{4}(n)\tag"(5.7)"\\
&\quad+{\frac {18}{11}}\,M_{6}(n)+\left (108\,{n}^{2}-24\,n+1\right )\,N_{2}(n),
\endalign
$$
$$
\align
N_{8}(n)&={\frac {2}{913}}\,\left (-72972\,{n}^{3}-1728\,{n}^{2}+5667\, n-289\right )\,M_{2}(n)\tag"(5.8)"\\
&+{\frac {280}{913}}\,\left (732\,{n}^{2}-195\,n+8 \right )\,M_{4}(n)+{\frac {84}{913}}\,\left (-196\,n+15\right )\,M_6(n)\\
&+{\frac {1248}{913}}\,M_{8}(n)+\left (-864\,{n}^{3}+360\,{n}^{2} -36\,n+1\right )\,N_{2}(n),
\endalign
$$
$$
\align
&N_{10}(n)\tag"(5.9)"\\
&={\frac {2}{5951847}}\,\left (3588144480\,{n}^{4}-805458600\, {n}^{3}-398007108\,{n}^{2}+56257647\,n-1794592\right )\,M_{2}(n)\\
&+{\frac  {140}{5951847}}\,\left (-72270360\,{n}^{3}+36826920\,{n}^{2}-3625245\,n+ 104002\right )\,M_{4}(n)\\
&+{\frac {210}{1983949}}\,\left (1421544\,{n}^{2} -380744\,n+13519\right )\,M_{6}(n)\\
&+{\frac {120}{1983949}}\,\left (- 282435\,n+18796\right )\,M_{8}(n)+{\frac {2724}{2173}}\,M_{10}(n) \\
&+\left (6480\,{n}^{4}-4320\,{n}^{3}+756\,{n}^{2}-48\,n+1\right )\,N_{2}( n) 
\endalign
$$
$$
\allowdisplaybreaks
\align
N_{14}(n)&={\frac {1}{4505033323132497}}\,\left (-655918847016750354240 \,{n}^{6}\right.\tag"(5.10)"\\
 & +584104439765983424400\,{n}^{5}-88193910587689930464\,{n}^{4} \\
 & - 51255985689606317364\,{n}^{3} +12889219681488512844\,{n}^{2}\\
 & \left. - 1033571808069319887\,n+23432656561492057\right )\,M_{2}(n)\\
 & +{\frac {364}{4505033323132497}}\,\left (2544016408481081520\,{n}^{5}\right.\\
 & -2986029950270749200\, {n}^{4}+1233083592931144500\,{n}^{3}\\
 & -185464100558325420\,{n}^{2}+ 12124758510318780\,n\\
 & \left. -229618708346923\right )\,M_{4}(n)\\
& +{\frac {728}{500559258125833}}\,\left (-12932704022040180\,{n}^{4}\right.\\
& +11781511098477120\,{n}^ {3}-3661921161131415\,{n}^{2}+234233352768436\,n\\
& \left. -7334109150929\right )\,M_{6}(n)+\\
& {\frac {364}{500559258125833}}\,\left (3327634333443960\,{n}^{3}\right.\\
& - 2184561928177200\,{n}^{2}+464283118670595\,n\\
& \left. -12774042869566\right )\,M_{8}(n)\\
& +{\frac {2002}{6030834435251}}\,\left (-758615153688\,{n}^{2}+ 404700708960\,n\right.\\
& \left. -24122003839\right )\,M_{10}(n)\\
& +{\frac {25388554464}{2775349487}}\,\left (n-1\right )\,M_{12}(n)\\
& +{\frac {139497552}{120667369 }}\,M_{14}(n)\\
&  +{\frac {1}{138}}\,\left (-107775360\,{n}^{6}+143700480\, {n}^{5}\right.\\
& \left. -70752528\,{n}^{4}+14978304\,{n}^{3}-1456488\,{n}^{2}+64320\,n-1045 \right )\,N_{2}(n)\\
& +{\frac {91}{138}}\,\left (-36\,n+13\right )\,N_{12}(n) 
\endalign
$$
\endproclaim
\demo{Proof} 
We first show that
relations such as \tagref{5.6}--\tagref{5.10} 
must exist. For $k\ge1$ we define
$$
\align
T_{k}
&= (2k-1)(k-1) R_{2k}
+ 6\sum_{i=1}^{k-1} \binom{2k}{2i}\left(2^{2i-1}-1\right) \delta_q(R_{2k-2i})\\
&\quad + \sum_{i=1}^{k-1}\left[
\binom{2k}{2i+2}\left(2^{2i+1}-1\right) - 2^{2i}\binom{2k}{2i+1} + \binom{2k}{2i}
\right] R_{2k-2i};
\tag"(5.11)"
\endalign
$$
i.e., $T_{k}$ is the function on the right side of \tagref{5.1} with $a=2k$.
By considering the left side of \tagref{5.1} and using Theorem \thmref{4.2} we
have
$$
T_{k} \in P\, \Cal W_{k};
\tag"(5.12)"
$$
i.e., there is a $\Phi\in\Cal W_{k}$ such that
$$
T_{k} = P \,\Phi.
\tag"(5.13)"
$$
For $k\ge1$ we consider the set of functions
$$
\Cal C_{k} = \{ \delta_q^m(C_{2j})\,:\, 1\le j\le k,\, j+m\le k\}
\subset P\,\Cal W_{k}.
\tag"(5.14)"
$$
Then 
$$
\left | \Cal C_{k} \right | = \frac{k(k+1)}{2}.
\tag"(5.15)"
$$
For $1\le k \le 5$ we observe from Corollary \thmref{3.6}
that
$$
\text{dim}\,\Cal W_k = \frac{k(k+1)}{2}.
\tag"(5.16)"
$$
Clearly,
$$
\dim{P\,\Cal W_k} = \dim{\Cal W_k},
\tag"(5.17)"
$$
for all $k$.
Thus, for $1 \le k \le 5$ there is a linear relation between
$T_{k}$ and the functions in $\Cal C_k$.
For $k=1$ the relation is $T_1=0$, which is trivial.
For each $2 \le k \le 5$ we have used \maple to find this relation.
This gives identities \tagref{5.6}--\tagref{5.9}.
Corollary \thmref{3.6} gives
$$
\dim \Cal W_6 = 22.
\tag"(5.18)"
$$
Since $\left|\Cal C_6\right|=21$ 
there is no reason to expect 
an analogous 
relation 
for $k=6$, and in fact 
calculation shows that the functions
in $\Cal C_6$ together with $T_{6}$ are linearly independent.
We note that
$$
T_{6},\, \delta_q(T_{6}),\, T_{7} \in P\,\Cal W_7.
\tag"(5.19)"
$$
Luckily we have     
$$
\dim \Cal W_7 = 30,\quad \text{and}\quad \dim \Cal C_7=28,
\tag"(5.20)"
$$
so there must exist a linear relation between
the three functions in \tagref{5.19} and the functions
in $\Cal C_7$. We used \maple to find this relation
and \tagref{5.10}.  \qed
\enddemo

Although \tagref{5.5} does not hold for $k=6$. There is a relation
for $k=6$ but involving an additional function. For an integer $r$
define $p_r(n)$ by
$$
\sum_{n\ge0} p_r(n) q^n = \prod_{n=1}^\infty (1-q^n)^r.
\tag"(5.21)"
$$
It is well-known that 
$$
\Delta(q) = \sum_{n\ge1} \tau(n) q^n =
\sum_{n\ge1} p_{24}(n-1) q^n = q\prod_{n=1}^\infty (1-q^n)^{24}
\tag"(5.22)"
$$
is a modular form of weight $12$. See \cite{S,pp.95--97}.
It follows that
$$
P \Delta = \sum_{n\ge1} p_{23}(n-1) q^n = q\prod_{n=1}^\infty (1-q^n)^{23}
\in P \Cal W_{6}.
\tag"(5.23)"
$$
By \tagref{5.18} and the fact that $\left|\Cal C_6\right|=21$ we see
that there must be a linear relation between $P\Delta$, $T_{12}$ and
the functions in $\Cal C_{6}$. We used \maple to find this
relation. It is given in the following 

\proclaim{Theorem 5.2}
For $n\ge1$ we have
$$
\align
&p_{23}(n-1)= \tag"(5.24)"\\
&
{\frac {1}{897196601564928}}\,\left (-57917897540518785552\,{n}^{5}+
30652078276547889552\,{n}^{4} \right.\\
&\qquad \left.+ 5952274737922797228\,{n}^{3}-
2214892612179680256\,{n}^{2}+188772676333745691\,n\right.\\
&\qquad \left. -4410708034409819
\right )\,M_{2}(n)
\\
&+
{\frac {5}{224299150391232}}\,\left (4089872889595634400\,{n}^{4}-
3320629034843596140\,{n}^{3}\right.\\
&\qquad \left. +593555423164294752\,{n}^{2}-
40741028214970311\,n+815166233039851\right )\,M_{4}(n)
\\
&+
{\frac {13}{7120607948928}}\,\left (-4612652508217680\,{n}^{3}+
2500384365901740\,{n}^{2}\right.\\
&\qquad \left. -190834728931028\,n+5730847932535\right )\,M_{6}(n)
\\
&+
{\frac {65}{24922127821248}}\,\left (431597256867684\,{n}^{2}-
112947999359631\,n+3472477850182\right )\,M_{8}(n)
\\
&+
{\frac {143}{600533200512}}\,\left (-555655003092\,n+33496841951
\right )\,M_{10}(n)+{\frac {16986177}{1919176}}\,M_{12}(n
)
\\
&+
{\frac {24599722121}{3316336128}}\,\left (-46656\,{n}^{5}+45360\,{n}^
{4}-12096\,{n}^{3}+1296\,{n}^{2}-60\,n+1\right )\,N_{2}(n)
\\
&
-{\frac {24599722121}{3316336128}}\,N_{12}(n).
\endalign
$$
\endproclaim

\vfill
\eject

\head 6. Congruence relations between rank and crank moments
\endhead

In the previous section we considered exact linear relations between
rank and crank moments. In this section we consider the analogous
problem modulo $p$.
One way to find such relations is to examine the denominators
of the rational numbers that occur in an exact relation.
For example, if we multiply both sides of \tagref{5.7} by $11$
and reduce modulo $11$ we find
$$
M_6(n) \equiv (2n+3) M_4(n) - (n+8)^2 M_2(n) \pmod{11}.
\tag"(6.1)"
$$
We could try the same sort of thing with \tagref{5.6}.
If we multiply both sides of \tagref{5.6} by $3$ and reduce modulo $3$
we obtain
$$
M_4(n) \equiv M_2(n) \pmod{3},
\tag"(6.2)"
$$
but this relation is trivial since $k^4\equiv k^2\pmod{3}$, for all integers
$k$.

There is another way in which a congruence between rank
and crank moments may arise. In \tagref{5.14} we defined
the set $\Cal C_k$. Let $k\ge1$. The set $\Cal C_k$
seems to be linearly independent over $\Bbb Q$. For small $k$ and 
certain primes $p$ the set $\Cal C_k$ is linearly dependent
over $\Bbb Z_p$. For example, consider the case $k=3$.
$$
\Cal C_3 =\{C_2,\delta_q(C_2),\delta_q^2(C_2),C_4,\delta_q(C_4),C_6\}.
$$
We want to find congruences between the elements of $\Cal C_3$.
We consider the $6\times6$ matrix $A$ whose $(i,j)$-th entry is
the coefficient of $q^i$ in the $q$-expansion of the $j$th element
of $\Cal C_3$.
$$
\align
A &= \pmatrix
M_{{2}}(1)& M_{{2}}(1)&M_{{2}}(1)&M_{{4}}(1)&M_{{4}}(1)&M_{{6}}(1)\\
M_{{2}}(2)&2\,M_{{2}}(2)&4\,M_{{2}}(2)&M_{{4}}(2)&2\,M_{{4}}(2)&M_{{6}}(2)\\
M_{{2}}(3)&3\,M_{{2}}(3)&9\,M_{{2}}(3)&M_{{4}}(3)&3\,M_{{4}}(3)&M_{{6}}(3)\\
M_{{2}}(4)&4\,M_{{2}}(4)&16\,M_{{2}}(4)&M_{{4}}(4)&4\,M_{{4}}(4)&M_{{6}}(4)\\
M_{{2}}(5)&5\,M_{{2}}(5)&25\,M_{{2}}(5)&M_{{4}}(5)&5\,M_{{4}}(5)&M_{{6}}(5)\\
M_{{2}}(6)&6\,M_{{2}}(6)&36\,M_{{2}}(6)&M_{{4}}(6)&6\,M_{{4}}(6)&M_{{6}}(6)
\endpmatrix \\
&=
\pmatrix
2&2&2&2&2&2\\
8&16&32&32&64&128\\
18&54&162&162&486&1458\\
40&160&640&544&2176&8320\\
70&350&1750&1414&7070&32710\\
132&792&4752&3300&19800&103092
\endpmatrix
\endalign
$$ 
We find that
$$
\det(A)=-110361968640 = -2^{17}\cdot3^7\cdot5\cdot7\cdot11.
$$
Since the determinant is nonzero the functions in $\Cal C_3$
are linearly independent. If there is a linear congruence 
for the functions in $\Cal C_3$ $\pmod{p}$, then $p$ will divide
this determinant. The occurrence of $p=11$ is confirmed by \tagref{6.1}.
We suspect that there may be a relation modulo $7$.
We find that
$$
(n+2) M_4(n) + (6n^2 + 4n + 1) M_2(n) \equiv 0 \pmod{7},
\tag"(6.3)"
$$
for all $n$. This is easily proved by writing the generating function
of the left side of \tagref{6.3} in terms of $P$, $\Phi_1$,
$\Phi_3$, $\Phi_5$:
$$
(\delta_q + 2) C_4 + (6\delta_q^2 + 4\delta_q + 1) C_2
=
-{\frac {7}{15}}\,{P}\,\left (180\,\Phi_{{1}}\Phi_{{3}}
-30\,{\Phi_{{1}}}^{2}-18\,\Phi_{{5}}-35\,\Phi_{{3}}-7\,\Phi_{{1}}\right )
\tag"(6.4)"
$$

We have used both methods described above to obtain congruences
between rank and crank moments. These are collected
together in the following
\proclaim{Theorem 6.1}
For $n\ge0$ we have
$$
(n+2) M_4(n) + (6n^2 + 4n + 1) M_2(n) \equiv 0 \pmod{7},
\tag"(6.5)"
$$
$$
\align
 (n + 5)^3 \,M_4(n) &\equiv
\left (5\,{n}^{4}+10\,{n}^{3}+8\,{n}^{2}+8\,n+9\right )\,M_{2}(n)
\pmod{11},
\tag"(6.6)"\\
M_6(n) &\equiv 2(n+7) M_4(n) - (n+8)^2 M_2(n) \pmod{11},
\tag"(6.7)"\\
M_{8}(n) &\equiv
2\,\left (n+5\right )\left ({n}^{2}+5\,n+10\right )\,M_{2}(n)+ 
6\left ({n}^{2}+n+1\right )\,M_{4}(n)
\pmod{11},
\tag"(6.8)"
\endalign
$$
$$
\align
M_{10}(n) &\equiv
4\,(n+7) \left ({n}^{3}+5\,{n}^{2}+29\,n+32\right )\,M_{2}(n)
\tag"(6.9)"\\
&\quad+39\,\left (n+7\right )\left (n+14\right )\left (n+39\right)\,M_{4}(n)+\left (6\,{n}^{2}+34\,n+39\right )\,M_{6}(n)\\
&\quad+35 \,\left (n+13\right )\,M_{8}(n)
\pmod{41},
\endalign
$$
$$
\align
& (n+7)\left  (n+25\right )\left (n+31\right )\left ({n}^{2}+28\,n+6\right )
\,N_{2}(n)+\,N_{12}(n) 
\tag"(6.10)"\\
&\equiv
\left (4\,{n}^{5}+10\,{n}^{4}+6\,{n}^{3}+30\,{n}^{2}+31\,n+33\right )\,M_{2}(n)
\\
&\quad+\left ({n}^{4}+4\,{n}^{3}+42\,{n}^{2}+24\,n+30\right )\,M_{4}(n)+40\,\left (n+10\right )\left ({n}^{2}+32\,n+7\right )\,M_{6}(n)\\
&\quad+31\,\left (n+13\right )\left (n+41\right )\,M_{8}( n)+\left (n+11\right )\,M_{10}(n)+22\,M_{12}(n)
\pmod{43},
\endalign
$$
$$
\align
M_{10}(n) &\equiv 36\,\left (n+19\right )\left ({n}^{3}+50\,{n}^{2}+20\,n+36\right )\, M_2(n)
\tag"(6.11)"\\
&\quad+\left (52\,{n}^{3}+28\,{n}^{2}+26\,n+52\right )\,M_{4}(n)+\left (36\,{n}^{2}+11\,n+32\right )\,M_{6}(n)\\
&\quad+47\, \left (n+17\right )\,M_{8}(n)\pmod{53},
\endalign
$$
$$
\align
M_8(n) &\equiv
\left (10\,{n}^{3}+73\,{n}^{2}+40\,n+82\right )\,M_{2}(n)+\left ( 72\,{n}^{2}+23\,n+28\right )\,M_{4}(n)
\tag"(6.12)"\\
&\quad+10\,\left (n+41\right )\,{M}_6(n)
\pmod{83},
\endalign
$$
$$
\align
&N_{12}(n) 
-367\,\left  (n+332\right )\left (n+487\right )\left (n+664\right )\left ({n}^{2}+ 265\,n+155\right )\,N_{2}(n) 
\tag"(6.13)"\\
&\equiv
352\,\left (n+247\right )\left (n+734\right )\left ({n}^{3}+147\,{n}^{2}+597\,n+363\right )\,M_{2}(n)\\
&\quad +88\,\left (n+530\right )\left (n+ 701\right )\left (n+709\right )\left (n+740\right )\,M_{4}(n)\\
&\quad+577 \,\left (n+114\right )\left (n+427\right )\left (n+682\right )\,M_{6}(n)+\left (295\,{n}^{2}+177\,n+674\right )\,M_{8}(n)\\
&\quad+271\, \left (n+336\right )\,M_{10}(n)+654\,M_{12}(n)
\pmod{797},
\endalign
$$
$$
\align
M_{14}(n)&\equiv
\left (44976165\,{n}^{6}+23584476\,{n}^{5}+19728425\,{n}^{4}+8711555\,
{n}^{3}+36781660\,{n}^{2}\right.
\tag"(6.14)"\\
&\qquad \left. +70780973\,n+108798274\right )\,M_{2}(n)
\\
&
+
77429163\,\left (n+4141548\right )\left (n+113894720\right )\left ({n}
^{3}+42853554\,{n}^{2}+28914352\,n\right.\\
&\qquad \left. +100598975\right )\,M_{4}(n)
\\
&
+
\left (12571854\,{n}^{4}+82951807\,{n}^{3}+9501843\,{n}^{2}+38248242\,
n+118847240\right )\,M_{6}(n)
\\
&
+
84218605\,\left (n+53645347 \right )\left ({n}^{2}+6688335\,n+93582728\right )
\,M_{8}(n)
\\
&
+
73449678\,\left (n+40889782\right )\left (n+59666357\right )\,M_{10}(n)
\\
&
+
89188917\,\left (n+120667368\right )\,M_{12}(n)
\\
&\pmod{120667369}.
\endalign
$$
\endproclaim
\demo{Proof}
We have already discussed \tagref{6.5} and \tagref{6.7}.
\tagref{6.12} follows from multiplying both sides of \tagref{5.8}
by $83$ and reducing modulo $83$. \tagref{6.6} and \tagref{6.11}
follow from multiplying both sides of \tagref{5.9} by
$11$ and $53$ respectively.
\tagref{6.14} follows from multiplying both sides of
\tagref{5.10} by $120667369$ and reducing modulo $120667369$.
The remaining relations follow by writing the 
generating functions for differences
between the left and right sides in terms of $P$, $\Phi_1$, $\Phi_3$,
and $\Phi_5$, and checking that the coefficient of each monomial
is divisible by the appropriate prime. \qed
\enddemo

We may multiply both sides of \tagref{5.10} by $23$ and then reduce
modulo $23$ to find a congruence for rank and crank moments.
A stronger relation follows from reducing \tagref{5.24}
modulo $23$ and using
$$
\sum_{n\ge1} p_{23}(n-1) q^n = q\prod_{n=1}^\infty (1-q^n)^{23}
\equiv 
 q\prod_{n=1}^\infty (1-q^{23n}) \pmod{23},
\tag"(6.15)"
$$
and Euler's pentagonal number theorem \cite{An, p.11}
$$
\prod_{n=1}^\infty (1-q^n) = \sum_{n=-\infty}^\infty (-1)^n q^{n(3n+1)/2}.
\tag"(6.16)"
$$
This gives
\proclaim{Theorem 6.2}
For $n\ge0$ we have
$$
\align
&
4\,\left ({n}^{2}+n+14\right )\left ({n}^{3}+{n}^{2}+15\right )\,{M}_2(n)
+
\left (10\,{n}^{4}+2\,{n}^{3}+8\,{n}^{2}+21\,n+22\right )\,M_{4}(
n)
\\
&
+
13\,\left (n+18\right )\left ({n}^{2}+21\,n+13\right )\,M_{6}(n)
+
5\,n\left (n+6\right )\,M_{8}(n)+15\,\left (n+19\right )\,{M}_{10}(n)
\\
&
+
M_{12}(n)
+
12\,\left (n+10\right )\left (n+14\right )\left (n+19\right )\left (n+
20\right )\left (n+21\right )\,N_{2}(n)
+
N_{12}(n)\\
& \equiv
\cases (-1)^k \pmod{23}, &\text{if $n=23k(3k\pm1)/2\,+\,1$},\\
0 \pmod{23}, &\text{otherwise}.
\endcases
\endalign
$$
\endproclaim

\subhead
Acknowledgement
\endsubhead
We would like to thank the referee for alerting us to the paper of
Kaneko and Zagier \cite{K-Z}.

\vfill
\eject

\enddocument